\documentclass[10pt]{amsart}
\usepackage{latexsym, amsmath}

 \renewcommand{\epsilon}{\varepsilon}
 \newcommand{\newsection}[1]
 {\subsection{#1}\setcounter{theorem}{0} \setcounter{equation}{0}
 \par\noindent}

 \newtheorem{theorem}{Theorem}
 
 \newtheorem{lemma}[theorem]{Lemma}
 \newtheorem{corr}[theorem]{Corollary}
 
 \newtheorem{proposition}[theorem]{Proposition}
 \newtheorem{deff}[theorem]{Definition}
 
 \newcommand{\bth}{\begin{theorem}}
 \newcommand{\ble}{\begin{lemma}}
 \newcommand{\bcor}{\begin{corr}}
 \newcommand{\bdeff}{\begin{deff}}
 \newcommand{\bprop}{\begin{proposition}}
 \newcommand{\eth}{\end{theorem}}
 \newcommand{\ele}{\end{lemma}}
 \newcommand{\ecor}{\end{corr}}
 \newcommand{\edeff}{\end{deff}}
 
 \newcommand{\eprop}{\end{proposition}}
  \newcommand{\cd}{\, \cdot\, }

 \renewcommand{\Pi}{\varPi}

 \renewcommand{\epsilon}{\varepsilon}

\renewcommand{\square}{\Box}

\begin{document}

 \title{Global existence for nonlinear wave equations with multiple speeds}
\thanks{The author was supported in part by the NSF}
\author{Christopher D. Sogge}
\address{Department of Mathematics, The Johns Hopkins University,
Baltimore, MD 21218}

 \maketitle

 \newsection{Introduction}

We shall be concerned with the Cauchy problem for quasilinear
systems in three space dimensions of the form
\begin{equation}\label{i.1}
\partial^2_tu^I-c^2_I\Delta u^I = C^{IJK}_{abc}\partial_c
u^J\partial_a\partial_b u^K + B^{IJK}_{ab}\partial_a u^J\partial_b
u^K, \quad I=1,\dots, D.
\end{equation}
 Here we are using the convention of summing repeated indices, and
$\partial u$ denotes the space-time gradient, $\partial
u=(\partial_0 u,
\partial_1 u, \partial_2 u, \partial_3u)$, with
$\partial_0=\partial_t$, and $\partial_j=\partial_{x_j}$,
$j=1,2,3$.   We shall be in the nonrelativistic case where we
assume that the wave speeds $c_k$ are all positive but not
necessarily equal.

The main difficulty in the nonrelativistic case is that one can
only use a smaller group of commuting vector fields.  In
particular, since one cannot use the generators of the hyperbolic
rotations, due to the different wave speeds, the earlier approach
of Klainerman \cite{K} breaks down.  This is because the invariant
Sobolev inequality that plays a key role in \cite{K} does not hold
if one uses a smaller collection of vector fields.  The conformal
approach of Christodoulou \cite{C} also does not seem to apply to
the nonrelativistic approach.  Also, unlike \cite{K}, our
techniques do not use Morawetz's conformal vector field.

In \cite{C1}, \cite{C}, \cite{K}, and \cite{K4} the
(3+1)-dimensional case was handled.  The null condition was first
identified and shown to lead to global existence of small
solutions.  Without the null condition, small solutions remain
smooth ``almost globally" \cite{JK}, but arbitrarily small
compactly supported initial data can develop singularities in
finite time \cite{J}.

We shall assume that the nonlinear terms satisfy a null condition.
Let us first assume, for simplicity, that the wave speeds $c_I$,
$I=1,\dots, D$ are distinct.  In this case, the null condition
only involves the self-interactions of each wave family. First we
require that self-interactions among the quasilinear terms satisfy
the standard null condition for the various wave-speeds:
\begin{equation}\label{i.2}
C^{IJJ}_{abc}\xi_a\xi_b\xi_c=0 \quad \text{whenever} \, \,
\frac{\xi_0^2}{c_J^2}-\xi_1^2-\xi_2^2-\xi_3^2=0, \, \,
I,J=1,\dots,D.
\end{equation}
We shall require that the self-interacting part of the semilinear
terms satisfy the standard null condition
\begin{equation}\label{i.3}
B^{IJJ}_{ab}\xi_a\xi_b=0, \quad \text{whenever} \, \,
\frac{\xi_0^2}{c_J^2}-\xi_1^2-\xi_2^2-\xi_3^2=0, \,  I,J=1,\dots,
D.
\end{equation}

If one allows repeated wave speeds, one must require that
interactions of families with the same speeds satisfy a null
condition.  Specifically, if we let $\mathcal{I}_p=\{I: \,
c_I=c_{I_p}, \, 1\le I\le D\}$, then the above null condition is
extended to be
\begin{equation}\label{i.2'}
C^{IJK}_{abc}\xi_a\xi_b\xi_c=0 \quad \text{whenever} \, \,
\frac{\xi_0^2}{c_{I_p}^2}-\xi_1^2-\xi_2^2-\xi_3^2=0, \, \,
(J,K)\in \mathcal{I}_p\times\mathcal{I}_p,\, 1\le I\le D,
\end{equation}
and
\begin{equation}\label{i.3'}
B^{IJK}_{ab}\xi_a\xi_b=0, \quad \text{whenever} \, \,
\frac{\xi_0^2}{c_{I_p}^2}-\xi_1^2-\xi_2^2-\xi_3^2=0, \, (J,K)\in
\mathcal{I}_p\times \mathcal{I}_p, \, 1\le I\le D.
\end{equation}

Since we are going to use the energy integral method, we also
require that the metric perturbation terms in the system
\eqref{i.1} are symmetric:
\begin{equation}\label{symmetry}
C^{IJK}_{abc}=C^{IJK}_{bac}=C^{KJI}_{abc}, \quad 0\le a,b,c\le 3,
\, \, 1\le I,J,K\le D.
\end{equation}

If the symmetry condition and these null conditions hold, we shall
show that \eqref{i.1} has a global solution, provided that the
initial data is small.

 To prove this result, we shall have to use estimates involving
various vector fields. We shall use the generators of
translations, $\partial$, the generators of Euclidean rotations
$\Omega=(\Omega_1,\Omega_2,\Omega_3)=x\wedge \nabla_x$, and the
scaling operator
\begin{equation}\label{i.4}
S=t\partial_t +r\partial_r = x^a\partial_a.
\end{equation}
Here $X=(x_0,x_1,x_2,x_3)=(t,x)$ denotes a point in
$\mathbb{R}^4$.  We shall denote these eight vector fields as
\begin{equation}\label{Gamma}
\Gamma=(\Gamma_0,\dots,\Gamma_7)=(\partial,\Omega,S),
\end{equation} and we shall use the multi-index notation
$$\Gamma^\alpha =\Gamma_{\alpha_m}\cdots\Gamma_{\alpha_1},$$
if $\alpha=(\alpha_1,\dots,\alpha_m)$, for a sequence of indices
$\alpha_i\in \{0,\dots,7\}$ of length $|\alpha|=m$.

The D'Alembertian will be the operator
\begin{equation}\label{i.5}
\square = \text{Diag}(\square^1,\dots,\square^D), \quad \text{with
} \, \square^I=\partial^2_t-c^2_I\Delta.
\end{equation}
We can then write \eqref{i.1} as
\begin{equation}\label{i.6}
\square u = N(\partial u, \partial^2u),
\end{equation}
where
\begin{equation}\label{i.7}
N^I(\partial u, \partial^2u)=C^{IJK}_{abc}\partial_c
u^J\partial_a\partial_b u^K + B^{IJK}_{ab}\partial_a u^J\partial_b
u^K, \quad I=1,\dots, D.
\end{equation}

To describe the solution space we let
$$\mathcal{H}^m(\mathbb{R}^3)=\{f\in L^2(\mathbb{R}^3; \mathbb{R}^D):
\, (\langle x\rangle\nabla_x)^\alpha f\in L^2, \quad |\alpha|\le
m\},$$ denote the weighted Sobolev space with norm
$$\|f\|_{\mathcal{H}^m(\mathbb{R}^3)}=\sum_{|\alpha|\le
m}\|(\langle x\rangle\nabla_x)^\alpha f\|_{L^2(\mathbb{R}^3)}.$$
Here, $\langle x\rangle=(1+|x|^2)^{1/2}$.

We can now state our global existence theorem for Minkowski space.

\begin{theorem}\label{theoremi.1}  Assume that the nonlinear terms
\eqref{i.7} satisfy \eqref{symmetry} as well as the null condition
\eqref{i.2'} and \eqref{i.3'}.  Then the initial value problem for
\eqref{i.6} with initial data
$$\partial u(0,\cd)\in \mathcal{H}^{m-1}(\mathbb{R}^3), \quad m\ge 10,$$
satisfying
\begin{equation}\label{i.8}\|\nabla_x
u(0,\cd)\|_{\mathcal{H}^{9}}+\|\partial_t
u(0,\cd)\|_{\mathcal{H}^9}<\varepsilon, \end{equation}
 with
$\varepsilon>0$ sufficiently small, has a unique global solution
satisfying $u(t,\cd)\in \mathcal{H}^{m}(\mathbb{R}^3)$ for every
$t>0$.
\end{theorem}

Sideris and Yu \cite{SY} proved the special case of this theorem
where the semilinear terms vanished identically.  Their approach
differed from much of the previous work since they did not use any
estimates arising from the fundamental solution of the
d'Alembertian.  A limitation of their approach, though, is that it
only leads to good pointwise control of second and higher
derivatives, and this explains why they were not able to handle
equations with semilinear terms.

Our approach is more in line with the original proof of Klainerman
\cite{K}.  A key difference, though, is that our main estimates
exploit the $1/\langle x\rangle$ decay of solutions of the wave
equation, as opposed to the $1/\langle t \rangle$ decay, which is
much more difficult to obtain.  We are able to exploit this weaker
decay because of a key pointwise estimate for solutions of the
inhomogeneous wave equation that is adapted to the $1/\langle
x\rangle$ decay and only uses the vector fields in \eqref{Gamma}.
This estimate was proved in an earlier joint paper with Keel and
Smith \cite{KSS3}, which proved almost global existence for
quasilinear equations using only the above vector fields in the
Minkowski space setting, as well as for the case of Minkowski
space minus star-shaped spatial obstacles. In a future paper, we
hope to show how the techniques from this paper can show that in
the obstacle setting one has global existence for the above
Dirichlet-wave equations when the null condition is satisfied.
This would extend results in our earlier joint paper with Keel and
Smith \cite{KSS1}, and give a proof that works directly in
Minkowski space, instead of relying on the conformal method.

This paper is organized as follows.  In the next section we shall
recall $L^2$ and pointwise estimates from \cite{KSS2} and
\cite{KSS3} that will play an important role in our arguments.
After that, we shall see consequences of the null condition and
prove some related auxiliary estimates that are adapted to the
$L^2$ and pointwise estimates that we shall use.  Then, in the
final section we shall prove Theorem \ref{theoremi.1}.

The author would like to thank Kunio Hidano for some very helpful
conversations, and Mark Keel for useful comments regarding an
early draft of this paper.

\newsection{Background: Pointwise estimates and  $L^2$ estimates}

As above
 $\{\Omega\}=\{\Omega_{ij}\}$,
 $\Omega_{ij}=x_i\partial_j-x_j\partial_i$, $1\le i<j\le 3$,
are the Euclidean ${\mathbb R}^3$ rotation operators, and the
scaling operator is
$S=t\partial_t+x\cdot\nabla_x=t\partial_t+r\partial_r$.
 Then we require the following  result from \cite{KSS3}.

 \begin{proposition}\label{mainprop}  If $w\in C^4$ and $(\partial_t^2-\Delta) w=F$ in
 $[0,t]\times {\mathbb R}^3$, and the Cauchy data of $w$ are $0$
 at $t=0$, then
 \begin{equation}\label{1}
 t|w(t,x)| \le C\int_0^t \int_{{\mathbb R}^3}\sum_{|\alpha|\le
3} |\Gamma^\alpha F(s,y)| \frac{dy ds}{1+|y|}.
 \end{equation}
 \end{proposition}

In \cite{KSS3} the following variant of of \eqref{1} was actually
proved:
$$
 t|w(t,x)| \le C\int_0^t \int_{{\mathbb R}^3}\sum_{|\alpha|\le
 2, j\le 1} |S^j \Omega^\alpha F(s,y)| \frac{dy ds}{|y|}.
$$
This estimate of course implies \eqref{1} when $F(s,y)=0$ for
$|y|<1$.  One gets the estimate for the case where $F(s,y)=0$ for
$|y|>2$ by applying the preceding case to a translation such as
$F(s,y_1,y_2,y_3+3)$.  (The translation introduces the constant
vector fields.)  Combining the two cases by a partition of unity
yields \eqref{1} in full generality.

In addition to the pointwise estimate \eqref{1}, we also require
an $L^2_tL^2_x$ estimate that is a simple consequences of the
energy inequality and Huygen's principle.


 \begin{proposition}\label{prop3.1} Suppose that $v$ solves the
 wave equation $(\partial_t^2-\Delta) v=G$ in ${\mathbb R}_+\times
{\mathbb R}^3$. Then there is a uniform constant $C$ so that
\begin{multline}\label{3.1}
\|(1+r)^{-1}v'\|_{L^2(\{(s,x): \, 0\le s\le t\})}
+\|(1+r)^{-2}v\|_{L^2(\{(s,x): \, 0\le s\le t\})}
\\
\le C\|v'(0,\cd)\|_2+C\int_0^t \|G(s,\cdot)\|_{L^2(\mathbb R^3)}\,
ds.
\end{multline}
\end{proposition}

Here, and in what follows, $v'$ denotes the space-time gradient of
$v$, i.e., $v'=(\partial_t v,\nabla_xv)$.

%

In addition to this $L^2_tL^2_x$ estimate we shall also of course
need the standard energy estimates for solutions of perturbed wave
equations
\begin{equation}\label{2.7}
(\partial_t^2-c_I^2\Delta) u^I+\sum_{K=1}^D\, \sum_{0\le j,k\le
3}\gamma^{IK,jk}\partial_j\partial_ku^K =F^I, \quad I=1,\dots,D,
\end{equation}
which satisfy the symmetry conditions
\begin{equation}\label{gammasymmetry}
\gamma^{IK,jk}=\gamma^{IK,kj}=\gamma^{KI,jk}, \quad 0\le j,k\le 3,
\, \, 1\le I,K\le D.
\end{equation}

The associated energy form then is $e_0=\sum_{I=1}^De^I,$ where
\begin{equation}\label{2.8}
e^I(u,t)=(\partial_0u^I)^2+\sum_{k=1}^3c^2_I(\partial_ku^I)^2+2\sum_{J=1}^D\sum_{k=0}^3\gamma^{IJ,0k}\partial_0u^I\partial_ku^J-\sum_{J=1}^D
\gamma^{IJ,jk}\partial_ju^I\partial_ku^J.
\end{equation}
 If we assume that
\begin{equation}\label{2.9}
\sum_{I,J=1}^D\sum_{j,k=0}^3 |\gamma^{IJ,jk}|\le
\frac12\min(1,c_I^2),
\end{equation}
it follows that
\begin{equation}\label{2.10}
\frac12\min(1,c_I^2) |\nabla_{t,x}u|^2 \le e(u).
\end{equation}
If we let $E(u,t)^2=\int_{\mathbb{R}^3}e(u,t)\, dx$ be the
associated energy, then, assuming \eqref{gammasymmetry} and
\eqref{2.9}, we also have the energy inequality
\begin{equation}\label{2.11}
\partial_t E(u,t)\le C\|F(t,\cd)\|_2+CE(u,t)\sum_{\substack{0\le j,k,l\le 3 \\ 1\le I,K\le D}}
\|\partial_l \gamma^{IK,jk}(t,\cd)\|_\infty,
\end{equation}
where $C$ is an absolute constant (only depending on the wave
speeds $c_I$).

If we use the following commutator relations $[\square, Z]=0$,
when $\{Z\}=\{\partial_j, \Omega_{ij}\}$, and $[\square,
S]=2\square$, where, as above, $S$ is the scaling vector field, we
see that this implies
\begin{multline}\label{2.12}
\sum_{|\alpha|\le M}\partial_t E(\Gamma^\alpha u,t)\le
C\sum_{|\alpha|\le M}\|\Gamma^\alpha F(t,\cd)\|_2
+\sum_{\substack{|\alpha|\le M \\ I,K,j,k}}\| \, [\Gamma^\alpha,
\gamma^{IK,jk}\partial_j\partial_k] u(t,\cd)\|_2
\\
+C\sum_{|\alpha|\le M}E(\Gamma^\alpha u,t)\sum_{\substack{0\le j,k,l\le 3 \\
1\le I,K\le D}}\|\partial_l \gamma^{IK,jk}(t,\cd)\|_\infty,
\end{multline}

\newsection{Null form bounds and auxiliary estimates}

Here we shall prove simple bounds for the null forms.  They must
involve the weight $\langle c_kt-r\rangle$ due to the fact that we
are not using the generators of Lorentz rotations.

\begin{lemma}\label{lemman.1}  Suppose that the nonlinear form
$N(\partial u,\partial^2u)$ satisfies the null condition
\eqref{i.2}-\eqref{i.3}.  Then
\begin{equation}\label{n.1}
|C^{IJJ}_{abc}\partial_c u\partial_a\partial_b v| \le C\langle
r\rangle^{-1}\bigl( |\Gamma u|\, |\partial^2 v| +|\partial u|\,
|\partial\Gamma v|\bigr) + C\frac{\langle c_Jt-r\rangle}{\langle
t+r\rangle} |\partial u|\, |\partial^2 v|.
\end{equation}
Also,
\begin{equation}\label{n.2}
|B^{IJJ}_{ab}\partial_a u\partial_b v| \le C\langle
r\rangle^{-1}\bigl(|\Gamma u|\, |\partial v| + |\partial u|\,
|\Gamma v|\bigr) + C\frac{\langle c_Jt-r\rangle}{\langle
t+r\rangle} |\partial u|\, |\partial v|.
\end{equation}
\end{lemma}

\noindent{\bf Proof of Lemma \ref{lemman.1}:}  Since
$\nabla_x=\frac{x}r\partial_r-\frac{x}{r^2}\wedge \Omega$.  So, if
we introduce the two operators $D^\pm=\frac12(\partial_t\pm
c_J\partial_r)$ and the null vectors $Y^\pm=(1,\pm x/c_Jr)$
associated to the wave speed $c_J$, we have
$$(\partial_t,\nabla_x)=(Y^-D^-+Y^+D^+)-(0,\frac{x}{c_Jr^2}\wedge
\Omega).$$ Since we can write
$$D^+=\frac{c_J}{c_Jt+r}S-\frac{c_Jt-r}{c_Jt+r}D^-,$$
the preceding formula can be rewritten as
$$\partial=Y^-D^--\frac{c_Jt-r}{c_Jt+r}Y^+D^-+\frac{c_J}{c_Jt+r}Y^+S-(0,\frac{x}{c_Jr^2}\wedge\Omega).$$
Consequently,
$$\partial =Y^-D^-+R,$$
where
$$|Ru|\le C\langle r\rangle^{-1}|\Gamma u|+\frac{\langle c_Jt-r\rangle}{t+r}|\partial u|.$$
Therefore, we have
\begin{multline*}
C^{IJJ}_{abc}\partial_cu\partial_a\partial_b v =
C^{IJJ}_{abc}\bigl( Y^-_{a}Y^-_{b}Y^-_{c}D^-u(D^-)^2v
\\
+R_cu\partial_{ab}v+Y^-_{c}D^-uR_a\partial_b v
+Y^{-}_{c}D^-uY^-_{a}D^-R_bv \bigr).
\end{multline*}
Note that $(Y^-_0)^2/c^2_J - (Y^-_1)^2-(Y^-_2)^2-(Y^-_3)^2=0$,
therefore by \eqref{i.2'} the first term in the right of the last
equation must vanish.  Therefore, the bounds for $R$ lead to
\eqref{n.1}.  The proof of \eqref{n.2} is similar.  \qed

Since we shall be proving estimates for scalar functions in the
rest of this section, let us abuse notation a bit by letting
$\square =
\partial^2_t-\Delta$ here.

\begin{lemma}\label{lemmaa.2}  If $h\in
C^\infty_0(\mathbb{R}_+\times \mathbb{R}^3)$ then
\begin{equation}\label{a.2}
\|\langle t-r\rangle\nabla h'(t,\cd)\|_2 \le C\sum_{|\alpha|\le
1}\|\Gamma^\alpha h'(t,\cd)\|_2 + C\|\langle t+r\rangle\square
h(t,\cd)\|_2.
\end{equation}
Also, if $0<\delta<1/2$ is fixed then
\begin{multline}\label{a.3}
\|h'(t,\cd)\|_{L^6(|x|\notin [(1-\delta)t,(1+\delta)t])}
\\
\le C\langle t\rangle^{-1}\Bigl(\sum_{|\alpha|\le
1}\|\Gamma^\alpha h'(t,\cd)\|_2 + C\|\langle t+r\rangle\square
h(t,\cd)\|_2\Bigr).
\end{multline}
\end{lemma}

\noindent{\bf Proof of Lemma \ref{lemmaa.2}:}  Inequality
\eqref{a.3} is a consequence of \eqref{a.2} since
$$\|h'(t,\cd)\|_{L^6(|x|\notin [(1-\delta)t,(1+\delta)t])}
\le C\|\nabla h'(t,\cd)\|_{L^2(|x|\notin
[(1-2\delta)t,(1+2\delta)t])}+C\langle
t\rangle^{-1}\|h'(t,\cd)\|_2.$$ Inequality \eqref{a.2} is
essentially in \cite{KS} (see Lemma 2.3 and Lemma 3.1 in
\cite{KS}).  The first step is to notice that one has the
elementary pointwise estimate
$$
\langle t-r\rangle\bigl( \, |\partial\partial_th(t,x)|+|\Delta
h(t,x)|\bigr) \le C\sum_{|\alpha|\le 1}|\partial \Gamma^\alpha
h(t,x)|+C\langle t+r\rangle|\square h(t,x)|,
$$
which leads to the $L^2$ bounds
$$\|\langle t-r\rangle(|\partial\partial_t h(t,\cd)|+|\Delta h(t,\cd)|)\|_2
\le C\sum_{|\alpha|\le 1}\|\Gamma^\alpha h'(t,\cd)\|_2 +
C\|\langle t+r\rangle\square h(t,\cd)\|_2.$$ Therefore, to finish,
we need to see that
\begin{equation}\label{aa.4}
\|\langle t-r\rangle\nabla_x^2 h(t,\cd)\|_2\le C\sum_{|\alpha|\le
1}\|\Gamma^\alpha h'(t,\cd)\|_2 + C\|\langle t+r\rangle\square
h(t,\cd)\|_2.
\end{equation}
To see this, we note that if we sum over $1\le i,j\le 3$ and
integrate by parts twice we have
\begin{multline*}
\int_{\mathbb{R}^3}\langle t-r\rangle^2
\sum_{i,j}|\partial_i\partial_j h|^2 dx
=\int_{\mathbb{R}^3}\langle t-r\rangle^2|\Delta h|^2 dx
\\
-\sum_{i,j}\Bigl(\int_{\mathbb{R}^3}\Bigl[(\partial_i\langle
t-r\rangle^2)(\partial_ju)\overline{(\partial_i\partial_j h )}
-(\partial_j\langle t-r\rangle^2)(\partial_j
h)\overline{(\partial^2_i h)}\Bigr] dx\Bigr).
\end{multline*}
The first term on the right has already been shown to be dominated
by the right side of \eqref{a.2}, and since an application of
Schwarz's inequality shows that the second term is dominated by
$$\Bigl(\int_{\mathbb{R}^3}\langle t-r\rangle^2 \sum_{i,j}|\partial_i\partial_j
h|^2 dx\Bigr)^{1/2} \, \|h'(t,\cd)\|_2,$$ we conclude that
\eqref{aa.4} must hold, which finishes the proof. \qed

The following result will be useful for dealing with waves
interacting at different speeds.

\begin{corr}\label{corra.3} Fix $c_1, c_2>0$ satisfying $c_1\ne c_2$.  Then
if $u,v\in C^\infty_0(\mathbb{R}\times\mathbb{R}^3)$,
\begin{align}\label{a.4}
\int_{\mathbb{R}^3}|&u''(t,x)|\, |v'(t,x)| \langle x\rangle^{-1}
dx
\\
&\le C\langle t\rangle^{-1}\bigl(\sum_{|\alpha|\le
1}\|\Gamma^\alpha u'(t,\cd)\|_2 + \|\langle
t+r\rangle(\partial_t^2-c_1^2\Delta) u(t,\cd)\|_2\bigr) \|\langle
x\rangle^{-1} v'(t,\cd)\|_2 \notag
\\
&+C\langle t\rangle^{-4/3}\bigl(\sum_{|\alpha|\le
1}\|\Gamma^\alpha u'(t,\cd)\|_2 + \|\langle
t+r\rangle(\partial_t^2-c_1^2\Delta) u(t,\cd)\|_2\bigr) \notag
\\
&\qquad\qquad\times
 \bigl(\sum_{|\alpha|\le
1}\|\Gamma^\alpha v'(t,\cd)\|_2+\|\langle
t+r\rangle(\partial_t^2-c_2^2\Delta) v(t,\cd)\|_2\bigr) \notag
\end{align}
\end{corr}

\noindent{\bf Proof:}  Let $2\delta <|c_1-c_2|$.  Then if we use
 Schwarz's inequality and \eqref{a.2} we see that we can bound
 $$\int_{|c_1t-r| >\delta t}
 |u''(t,x)|\,
 |v'(t,x)| \langle x\rangle^{-1} dx$$
 by the first term in the right side of \eqref{a.4}.

For the next step we split the remaining region where $|r-c_1t|\le
\delta t$ into the annuli, $A_j$, where $\langle c_1t-r\rangle\in
[2^j,2^{j+1})$, $j=0,1,2,\dots$.  Assuming that $A_j\cap\{x:
|r-c_1t|\le \delta t\}\ne \emptyset$, we can use H\"older's
inequality to find that
\begin{multline*}\int_{\langle c_1t-r\rangle\in [2^j,2^{j+1})} |u''(t,x)| \, |v'(t,x)|
\langle x\rangle^{-1} dx
\\
\le C t^{-1/3}2^{j/3}\|u''(t,\cd)\|_{L^2(\langle c_1t-r\rangle\in
(2^j, 2^{j+1}))}\|v'(t,\cd)\|_{L^6(\langle c_1t-r\rangle\in (2^j,
2^{j+1}))}.
\end{multline*}
Since $2\delta <|c_1-c_2|$, on the set where $|r-c_1t|\le \delta
t$ we have the lower bound $|r-c_2t|\ge \delta t$.
Therefore, we can apply \eqref{a.2} and \eqref{a.3} to see that
the right side is bounded by $2^{-2j/3}$ times the second term in
the right side of \eqref{a.4}, which, after summing over $j$,
implies that when we restrict the integration in the left side of
\eqref{a.4} to the the set where $|r-c_1t|\le \delta t$
the resulting expression is dominated by the second term in the
right of \eqref{a.4}.  Therefore, after summing over $j$, we also
have control of the analog of (3.6) where the integration is over
the region where $|c_1t-r|\le \delta$, which completes the proof.
\qed

To handle same-speed interactions, we shall need the following
similar result.

\begin{corr}\label{corra.4} Let $u,v\in C^\infty_0(\mathbb{R}\times\mathbb{R}^3)$.
 Then,
\begin{align}\label{a.5}
\int_{\mathbb{R}^3}\langle x\rangle^{-2}&|\partial^2 u(t,x)|\, |
v(t,x)|\, dx \\
&\le C\langle
t\rangle^{-1}\bigl(\sum_{|\alpha|\le1}\|\Gamma^\alpha
u'(t,\cd)\|_2+\|\langle t+r\rangle\square u(t,\cd)\|_2\bigr)
\|\langle x\rangle^{-2} v\|_2 \notag
\\
&+C\langle t\rangle^{-4/3}\bigl(\sum_{|\alpha|\le1}\|\Gamma^\alpha
u'(t,\cd)\|_2+\|\langle t+r\rangle\square
u(t,\cd)\|_2\bigr)\|v'\|_2. \notag
\end{align}
\begin{align}\label{aa.5}
&\int_{\mathbb{R}^3}\langle x\rangle^{-2}|u'(t,x)|\, |
v'(t,x)|\, dx \\
&\le C\log(2+t)\langle
t\rangle^{-1}\bigl(\sum_{|\alpha|\le1}\|\Gamma^\alpha
u'(t,\cd)\|_2+\|\langle t+r\rangle\square u(t,\cd)\|_2\bigr)
\|\langle x\rangle^{-1} v'\|_2 \notag
\\
&+C\langle t\rangle^{-2}\|u'(t,\cd)\|_2\|v'(t,\cd)\|_2. \notag
\end{align}
\begin{multline}\label{a.6}
\int_{\mathbb{R}^3}\frac{\langle t-r\rangle}{\langle
t+r\rangle}|\partial^2 u(t,x)|\, |v'(t,x)|\, \langle x\rangle^{-1}
\, dx
\\
 \le
C\langle t\rangle^{-1}\bigl(\sum_{|\alpha|\le1}\|\Gamma^\alpha
u'(t,\cd)\|_2+\|\langle t+r\rangle\square
u(t,\cd)\|_2\bigr)\|\langle x\rangle^{-1}v'(t,\cd)\|_2.
\end{multline}
\end{corr}

\noindent{\bf Proof of Corollary \ref{corra.4}:}  To prove
\eqref{a.5} we first notice that by using Schwarz's inequality and
(3.3) we get that
\begin{multline*}\int_{|x|\notin
[t/2,3t/2]}\langle x\rangle^{-2}|\partial^2 u(t,x)|\, | v(t,x)|\,
dx
\\
\le C\langle t\rangle^{-1}\bigl(\sum_{|\alpha|\le1}\|\Gamma^\alpha
u'(t,\cd)\|_2+\|\langle t+r\rangle\square u(t,\cd)\|_2\bigr)
\|\langle x\rangle^{-2} v\|_2.
\end{multline*}
If we use \eqref{a.2} and argue as in the proof of Corollary
\ref{corra.3} we can estimate the integral over $|x|\in
[t/2,3t/2]$. We do so by noting that for $j=0,1,2,\dots$ the
integral over $\{x: \, |x|\in [t/2,3t/2]\} \cap \{x: \, \langle
t-r\rangle \in [2^j, 2^{j+1})\}$ is dominated by
$$
\langle t\rangle^{-2}2^{-j}\bigl(\sum_{|\alpha|\le
2}\|\Gamma^\alpha u'(t,\cd)\|_2+\|\langle t+r\rangle\square
u(t,\cd)\|_2\bigr) t^{2/3}2^{j/3}\|v(t,\cd)\|_6.$$
We conclude
that the integral over the region where $|x|\in [t/2, 3t/2]$  is
dominated by the other term in the right side of \eqref{a.5} after
summing over $j$ and applying Sobolev's theorem.

To prove \eqref{aa.5} we note that it suffices to show that
\begin{multline*}
\int_{|x|<(t+1)/2}\langle x\rangle^{-2}|u'(t,x)|\, | v'(t,x)|\, dx
\\
\le C\log(2+t)\langle
t\rangle^{-1}\bigl(\sum_{|\alpha|\le1}\|\Gamma^\alpha
u'(t,\cd)\|_2+\|\langle t+r\rangle\square u(t,\cd)\|_2\bigr)
\|\langle x\rangle^{-1} v'\|_2,
\end{multline*}
since clearly the integral over $|x|>(t+1)/2$ is dominated by the
second term in the right side of \eqref{aa.5}.  However, if we use
H\"older's inequality and (3.4), we see that the integral over
$|x|<(t+1)/2$ is dominated by
\begin{multline*}
\log(2+t)\|u'(t,\cd)\|_{L^6(|x|<(t+1)/2)}\|\langle
x\rangle^{-1}v'\|_2
\\
\le C\log(2+t)\langle t\rangle^{-1}\bigl(\sum_{|\alpha|\le
1}\|\Gamma^\alpha u'(t,\cd)\|_2+\|\langle t+r\rangle\square
u(t,\cd)\|_2\bigr)\|\langle x\rangle^{-1}v'\|_2,
\end{multline*}
which handles the remaining part of \eqref{aa.5}.

To prove \eqref{a.6} we just use   Schwarz's inequality and (3.3)
to see that its left side is dominated by
\begin{multline*}
\langle t\rangle^{-1}\|\langle
t-r\rangle\partial^2u(t,\cd)\|_2\|\langle
x\rangle^{-1}v'(t,\cd)\|_2
\\
\le C\langle t\rangle^{-1}\bigl(\sum_{|\alpha|\le1}\|\Gamma^\alpha
u'(t,\cd)\|_2+\|\langle t+r\rangle\square
u(t,\cd)\|_2\bigr)\|\langle x\rangle^{-1}v'(t,\cd)\|_2,
\end{multline*}
 which completes the
proof. \qed

We shall also need the following simple result.

\begin{lemma}\label{lemma3.5} Let $u,v\in C^\infty_0(\mathbb{R}\times\mathbb{R}^3)$.
 Then,
\begin{multline}\label{3.10}
\|\, \langle t+r\rangle u(t,\cd)v(t,\cd)\|_{L^2(\mathbb{R}^3)} \le
C\|u(t,\cd)\|_{L^2(\mathbb{R}^3)}\, \langle
t\rangle\|v(t,\cd)\|_{L^\infty(\mathbb{R}^3)}
\\
+C\|u(t,\cd)\|_{L^2(\mathbb{R}^3)}\sum_{|\alpha|+|\beta|\le
2}\|\Omega^\alpha \partial_x^\beta v(t,\cd)\|_{L^2(\mathbb{R}^3)}.
\end{multline}
\end{lemma}

The proof is simple.  We first note that if we replace $\langle
t+r\rangle$ by $\langle t\rangle$ in the left side of \eqref{3.10}
the resulting quantity is dominated by the first term in the right
side of \eqref{3.10}. So to finish the proof, it would suffice to
show that
$$\|\, r u(t,\cd)v(t,\cd)\|_{L^2(\mathbb{R}^3)} \le
C\|u(t,\cd)\|_{L^2(\mathbb{R}^3)}\sum_{|\alpha|+|\beta|\le
2}\|\Omega^\alpha \partial_x^\beta
v(t,\cd)\|_{L^2(\mathbb{R}^3)}.$$ Since $|v(t,x)|\le
C\sum_{|\beta|\le 2}\|\partial^\beta_x v(t,\cd)\|_2$, we can prove
the analog of this estimate where the norm in the left is taken
over $|x|\le 2$. The $L^2$ bound over the region $|x|>2$ follows
from the following Sobolev estimate for $\mathbb{R}_+\times S^2$
$$|v(t,r\omega)|\le C\sum_{|\alpha|+j\le
2}\Bigl(\int_{|r-\rho|\le 1} \int_{S^2} |\Omega^\alpha
\partial_\rho^j v(t,\rho\omega)|^2 \, d\omega d\rho\Bigr)^{1/2},$$
since, when expressed in polar coordinates, the standard volume
form on $\mathbb{R}^3$ is $4\pi \rho^2d\rho d\omega$.

\newsection{Proof of Theorem 1.1}


We shall assume that the data satisfies the smallness condition
(1.12).  We then  wish to show that \eqref{i.6} has a global
solution if the null condition \eqref{i.2'}-\eqref{i.3'} holds.
For simplicity, we shall assume that the wave speeds $c_I$,
$I=1,\dots,D$, are distinct.  A simple modification of the
arguments to follow handles the general case where this assumption
is removed.

To proceed, we shall need to use a standard local existence
theorem:

\begin{theorem}\label{theorem4.1}  Suppose that $f$ and $g$ as
above
satisfy \eqref{i.8}.
Suppose also that the symmetry condition \eqref{symmetry} holds.
Then there is a $T>0$ so that the initial value problem
\eqref{i.6} with initial data has a $C^2$ solution satisfying
\begin{equation}\label{4.2}
u\in L^\infty([0,T]; H^{10}(\mathbb{R}^3)) \cap C^{0,1}([0,T];
H^9(\mathbb{R}^3)).
\end{equation}
The supremum of all such $T$ is equal to the supremum of all $T$
such that the initial value problem has a $C^2$ solution with
$\partial^\alpha u$ bounded for $|\alpha|\le 2$.
\end{theorem}

This result is essentially Theorem 6.4.11 in \cite{H}.  The latter
result, though, is just for scalar wave equations; however, the
same proof, which is based on energy inequalities, yields Theorem
\ref{theorem4.1} since we are assuming that the symmetry condition
\eqref{symmetry} and thus one can use the energy-integral method
exactly as in \cite{H}.

\bigskip

We now turn to the proof of Theorem 1.1.  We let $\varepsilon>0$
be as in \eqref{i.8}, and assume that we already have a $C^2$
solution of our equation for $0\le t\le T$ such that for such $t$
and small $\varepsilon$
\begin{equation}\label{4.3}
(1+|t|)\sum_{|\alpha|\le 4}|\Gamma^\alpha u'(t,x)|\le
A_0\varepsilon
\end{equation}
\begin{multline}\label{4.4}
\sum_{|\alpha|\le 9}\|\Gamma^\alpha
u'(t,\cd)\|_{L^2(\mathbb{R}^3)}+\sum_{|\alpha|\le 8}\|\langle
x\rangle^{-1}\Gamma^\alpha u'\|_{L^2(\{(s,x)\in [0,t]\times
\mathbb{R}^3\})}
\\
+\sum_{|\alpha|\le 8}\|\langle x\rangle^{-2}\Gamma^\alpha
u\|_{L^2(\{(s,x)\in [0,t]\times \mathbb{R}^3\})}
 \le
A_1\varepsilon(1+t)^{A_2\varepsilon}\sum_{|\alpha|\le
9}\|\Gamma^\alpha u'(0,\cd)\|_2.
\end{multline}
Clearly both estimates are valid if $T$ is small.

We then let $A_0$ be so large that \eqref{4.3} holds with $A_0$
replaced by $A_0/3$ if $u$ is replaced by the solution of the wave
equation $\square u_0=0$ with Cauchy data $(f,g)$, and $\square$
is as in \eqref{i.5}.  We shall then prove for $\varepsilon$
smaller than some number depending on $A_1$ and $A_2$ that

{\bf i)} \eqref{4.3} is valid with $A_0$ replaced by $A_0/2$;

{\bf ii)} \eqref{4.4} is a consequence of \eqref{4.3} for suitable
$A_1$, $A_2$.
\newline
By the local existence theorem it will follow that a solution
exists for all $t\ge 0$ if $\varepsilon$ is small enough.

\bigskip
\noindent{\bf Proof of i):}  Since the Cauchy data of
$\Gamma^\alpha u-\Gamma^\alpha u_0$ is $O(\varepsilon^2)$, it
suffices by Proposition \ref{mainprop} to prove that for small
$\varepsilon$
$$\sum_{|\alpha|\le 4}\sum_{|\beta|\le 3}
\iint_{0<s<T}|\Gamma^\beta \square \Gamma^\alpha
u'(s,y)|\frac{dyds}{1+|y|}\le C\varepsilon^2.$$ Using the
commutativity relations of the $\Gamma$ and $\square$, we can
write $\Gamma^\alpha \square \Gamma^\alpha u'$ as a sum of terms
of the form $\Gamma^\sigma \square u'$ with $|\sigma|\le
|\alpha|+|\beta|\le 7$.  Therefore, it suffices to prove that
\begin{equation}\label{4.5}
\sum_{|\alpha|\le 7}\iint_{0<s<T}|\Gamma^\alpha \square
 u'(s,y)|\frac{dyds}{1+|y|}\le C\varepsilon^2.
 \end{equation}
 To do so we note that the $I$-th component of $\Gamma^\alpha \square
 u'$, $|\alpha|\le 7$ is a linear combination of terms of the form
\begin{equation}\label{4.6}\sum_{a,b,c=0}^3\Gamma^\alpha
\partial\bigl(C^{IJK}_{abc}\partial_c u^J\partial_a\partial_b
u^K\bigr) + \sum_{a,b=0}^3\Gamma^\alpha
\partial\bigl(B^{IJK}_{ab}\partial_a u^J\partial_b u^K\bigr),
\quad |\alpha|\le 7.
\end{equation}

Let us first handle the contribution to \eqref{4.5} of the terms
in \eqref{4.6} with $J\ne K$.  If we use Corollary \ref{corra.3},
we find that if we fix $s\in (0,T)$ then for $J\ne K$ and
$|\alpha|\le 7$ we have
\begin{align}
\int_{\mathbb{R}^3}&\bigl(
|\Gamma^\alpha\partial(C^{IJK}_{abc}\partial_cu^J\partial_a\partial_bu^K)(s,y)|
+|\Gamma^\alpha
\partial(B^{IJK}_{ab}\partial_a u^J\partial_b u^K)|
\bigr)\frac{dy}{1+|y|}\label{4.7}
\\
&\le C \langle s\rangle^{-1}\bigl(\sum_{|\alpha|\le
9}\|\Gamma^\alpha u'(s,\cd)\|_2+\sum_{|\alpha|\le 8}\|\langle
s+r\rangle(\partial^2_t-c^2_J\Delta)\Gamma^\alpha
u^J(s,\cd)\|_2\bigr)\notag
\\ &\qquad \times \sum_{|\alpha|\le 8}  \|\langle r\rangle^{-1}\Gamma^\alpha u'(s,\cd)\|_2
\notag
\\
&+C\langle s\rangle^{-4/3}\bigl(\sum_{|\alpha|\le
9}\|\Gamma^\alpha u'(s,\cd)\|_2+\sum_{|\alpha|\le 8}\|\langle
s+r\rangle(\partial^2_t-c^2_J\Delta)\Gamma^\alpha
u^J(s,\cd)\|_2\bigr) \notag
\\ &\qquad \times\bigl(\sum_{|\alpha|\le 9}\|\Gamma^\alpha
u'(s,\cd)\|_2+\sum_{|\alpha|\le 8}\|\langle
s+r\rangle(\partial^2_t-c^2_K\Delta)\Gamma^\alpha
u^K(s,\cd)\|_2\bigr). \notag
\end{align}
Note that
$$|(\partial^2_t-c^2_J\Delta)\Gamma^\alpha u^J(s,y)|\le
C\sum_{|\beta|\le 4}|\Gamma^\beta u'(s,y)|\sum_{|\beta|\le
9}|\Gamma^\beta u'(s,y)|, \quad |\alpha|\le 8.$$ Therefore, if we
use Lemma \ref{lemma3.5} and our assumptions
\eqref{4.3}-\eqref{4.4} we conclude that
$$\sum_{|\alpha|\le
8}\|\langle s+r\rangle(\partial^2_t-c^2_J\Delta)\Gamma^\alpha
u^J(s,\cd)\|_2 \le C\varepsilon^2 \langle
s\rangle^{2A_2\varepsilon}.$$ Because of this, if we use
\eqref{4.3}-\eqref{4.4} again, we find that the left side of
\eqref{4.7} is dominated by
$$C\varepsilon\langle s\rangle^{-1+2A_2\varepsilon} \sum_{|\alpha|\le 8}  \|\langle x\rangle^{-1}u'(s,\cd)\|_2
+C\varepsilon^2\langle s\rangle^{-4/3+2A_2\varepsilon},$$ which
means that when $J\ne K$
\begin{align}\label{4.8}
\sum_{|\alpha|\le 7}&\iint_{0<s<T}\bigl(
|\Gamma^\alpha\partial(C^{IJK}_{abc}\partial_cu^J\partial_a\partial_bu^K)(s,y)|
+|\Gamma^\alpha
\partial(B^{IJK}_{ab}\partial_a u^J\partial_b u^K)|
\bigr)\frac{dyds}{1+|y|} \\
&\le C\int_0^T \bigl(\varepsilon\langle
s\rangle^{-1+2A_2\varepsilon} \sum_{|\alpha|\le 8}  \|\langle
r\rangle^{-1}u'(s,\cd)\|_2 +\varepsilon^2\langle
s\rangle^{-4/3+2A_2\varepsilon}\bigr) \, ds\notag
\\
&\le C_\delta\varepsilon\sum_{|\alpha|\le 8} \|\langle
s\rangle^{-\frac12+2A_2\varepsilon+\delta}\langle
r\rangle^{-1}\Gamma^\alpha u'\|_{L^2(\{(s,x)\in [0,t]\times
\mathbb{R}^3\})} +C\varepsilon^2 , \quad \delta>0,\notag
\end{align}
using the Schwarz inequality in the last step.  If we assume that
$\varepsilon$ and $\delta$ are small enough $\frac12
-3A_2\varepsilon-\delta>0$ then we can bound the first term in the
right by $C\varepsilon^2$ as well by using \eqref{4.4}.  One first
controls the $L^2$ norms when $s$ lies in a dyadic subinterval of
$[0,t]$ and then sums.

We therefore conclude that the terms in \eqref{4.6} with $J\ne K$
satisfy the bounds in \eqref{4.5}.  So to finish we have to
consider the terms with $J=K$, in which case we shall need to use
the null condition.

We first estimate the contribution of the quasilinear terms
satisfying the null condition.  If we use Lemma \ref{lemman.1} we
conclude that
\begin{align}\label{4.9}
\sum_{|\alpha|\le 7}\Bigl(|\Gamma^\alpha
\partial\bigl(C^{IJJ}_{abc}\partial_c&u^J\partial_a\partial_b u^J\bigr)|
+|\Gamma^\alpha
\partial\bigl(B^{IJJ}_{ab}\partial_au^J\partial_b u^J\bigr)|\Bigr)
\\
 &\le C\langle r\rangle^{-1} \Bigl[
\sum_{|\alpha|\le 8}|\Gamma^\alpha u^J| \sum_{|\alpha|\le
8}|\partial^2 \Gamma^\alpha u^J|+\sum_{|\alpha|\le 8}|\partial
\Gamma^\alpha u^J|^2 \Bigr] \notag
\\
&+C\frac{\langle c_Js-r\rangle}{\langle s+r\rangle}
\sum_{|\alpha|\le 8}|\partial\Gamma^\alpha u^J|\sum_{|\alpha|\le
8}|\partial^2 \Gamma^\alpha u^J|. \notag
\end{align}
If we use \eqref{a.5} we can handle the first term on the right.
Indeed,

\begin{align}\label{4.10}
\int_{\mathbb{R}^3}& \sum_{|\alpha|\le 8}|\partial^2 \Gamma^\alpha
u^J(s,y)|\sum_{|\alpha|\le 8}|\Gamma^\alpha u^J(s,y)|
\frac{dy}{(1+|y|)^2}
\\
&\le C\langle s\rangle^{-1}\Bigl(\sum_{|\alpha|\le
9}\|\Gamma^\alpha u'(s,\cd)\|_2 +\sum_{|\alpha|\le 8}\|\langle
s+r\rangle(\partial^2_t-c_J^2\Delta)\Gamma^\alpha
u^J(s,\cd)\|_2\Bigr) \notag
\\
&\qquad\qquad\times\sum_{|\alpha|\le 8}\|\langle
r\rangle^{-2}\Gamma^\alpha u(s,\cd)\|_2 \notag
\\
&+C\langle s\rangle^{-4/3}\Bigl(\sum_{|\alpha|\le
9}\|\Gamma^\alpha u'(s,\cd)\|_2 +\sum_{|\alpha|\le 8}\|\langle
s+r\rangle(\partial^2_t-c_J^2\Delta)\Gamma^\alpha
u^J(s,\cd)\|_2\Bigr) \notag
\\
&\qquad\qquad\times\sum_{|\alpha|\le 8}\|\Gamma^\alpha
u'(s,\cd)\|_2.  \notag
\end{align}
Using \eqref{aa.5} we can bound the second term in the right side
of \eqref{4.9}:
\begin{align}\label{4.11}
\int_{\mathbb{R}^3}&\sum_{|\alpha|\le 8}|\partial \Gamma^\alpha
u^J(s,y)|^2 \frac{dy}{(1+|y|)^2}
\\
&\le C\log(2+s)\langle s\rangle^{-1}\Bigl(\sum_{|\alpha|\le
9}\|\Gamma^\alpha u'(s,\cd)\|_2 +\sum_{|\alpha|\le 8}\|\langle
s+r\rangle(\partial^2_t-c_J^2\Delta)\Gamma^\alpha
u^J(s,\cd)\|_2\Bigr) \notag
\\
&\qquad\qquad\times\sum_{|\alpha|\le 8}\|\langle
r\rangle^{-1}\Gamma^\alpha u'(s,\cd)\|_2 \notag
\\
&+C\langle s\rangle^{-2}\sum_{|\alpha|\le 8}\|\Gamma^\alpha
u'(s,\cd)\|_2^2.\notag
\end{align}
And, similarly, if we use \eqref{a.6} we can estimate the last
term in \eqref{4.9}:
\begin{align}\label{4.12}
\int_{\mathbb{R}^3}&\frac{\langle c_Js-r\rangle}{\langle
s+r\rangle}\ \sum_{|\alpha|\le 8}|\partial^2 \Gamma^\alpha
u^J(s,y)|\, \sum_{|\alpha|\le 8}|\partial\Gamma^\alpha u^J(s,y)|
 \frac{dy}{1+|y|}
 \\
&\le C\langle s\rangle^{-1}\bigl(\sum_{|\alpha|\le
9}\|\Gamma^\alpha u'(s,\cd)\|_2+\sum_{|\alpha|\le 8} \|\langle
s+r\rangle(\partial^2_t-c_J^2\Delta)\Gamma^\alpha
u^J(s,\cd)\|_2\bigr) \notag
\\
&\qquad \qquad \times \sum_{|\alpha|\le 8}\|\langle
r\rangle^{-1}\Gamma^\alpha u'(s,\cd)\|_2. \notag
\end{align}

If we combine \eqref{4.9}-\eqref{4.12}, then the arguments used to
prove \eqref{4.8} yield
\begin{equation*}\sum_{|\alpha|\le
7}\iint_{0<s<T}
\Bigl(|\Gamma^\alpha\partial(C^{IJJ}_{abc}\partial_cu^J\partial_a\partial_bu^J)|
+|\Gamma^\alpha
\partial(B^{JJ}\partial_au^J\partial_b u^J)|
 \Bigr)\frac{dyds}{1+|y|}
 \le C\varepsilon^2,
 \end{equation*}
  as desired.  This, along with \eqref{4.8} yields \eqref{4.5}.
  Therefore, we have completed step i) of the proof.

\bigskip
\noindent{\bf Proof of ii):}  We need to apply \eqref{2.12} with
$M=9$, and
\begin{align*}
F^I &= \sum_{J,K}\sum_{a,b}B^{IJK}_{ab}\partial_au^J\partial_bu^K
\\
\gamma^{IK,ab}&=-\sum_J\sum_c C^{IJK}_{abc}\partial_c u^J.
\end{align*}
We then have \eqref{2.9}, assuming that $\varepsilon$ in
\eqref{4.3} is small.
 Since
$$\sum_{|\alpha|\le 9}\Bigl(|\Gamma^\alpha F|
+|[\Gamma^\alpha, \gamma^{IK,ab}\partial_a\partial_b]u|\Bigr) \le
\sum_{|\alpha|\le 4}|\Gamma^\alpha u'|\sum_{|\alpha|\le
9}|\Gamma^\alpha u'|,$$ we conclude from \eqref{4.3} that the
first two terms in the right side of \eqref{2.12} are
$$\le C\varepsilon(1+t)^{-1}\sum_{|\alpha|\le 9}E(\Gamma^\alpha
u,t).$$ Since \eqref{4.3} also implies that this must be the case
for the last term in \eqref{2.12}, we conclude that \eqref{4.3}
implies that
$$\sum_{|\alpha|\le 9}\partial_t E(\Gamma^\alpha u,t) \le
C\varepsilon(1+t)^{-1} \sum_{|\alpha|\le 9}E(\Gamma^\alpha u,t),$$
and hence that
$$\sum_{|\alpha|\le 9}E(\Gamma^\alpha u,t) \le
(1+t)^{C\varepsilon} \sum_{|\alpha|\le 9}E(\Gamma^\alpha u,0).$$
This implies that the first term in \eqref{4.4} satisfies the
desired bounds.

If we take
$$G^I=\sum_{J,K}\sum_{a,b}B^{IJK}_{ab}\partial_au^J\partial_b u^K
+\sum_{J,K}\sum_{a,b,c}C^{IJK}_{abc}\partial_c
u^J\partial_{ab}u^K,$$ then we can apply Proposition \ref{prop3.1}
to bound the other two terms in the left side of \eqref{4.4}.
Indeed, they  are controlled by
$$\sum_{|\alpha|\le 8}\|\Gamma^\alpha u'(0,\cd)\|_2
+\int_0^t \sum_{|\alpha|\le 8}\|\Gamma^\alpha G(s,\cd)\|_2\, ds.$$
We can control the last term if we use \eqref{4.3} and the fact
that we have shown that the first term in \eqref{4.4} is under
control.  By doing this we get
\begin{equation*}\sum_{|\alpha|\le 8}\|\Gamma^\alpha G(s,\cd)\|_2
\le C\sum_{|\alpha|\le 4}\|\Gamma^\alpha u'(s,\cd)\|_\infty
\sum_{|\alpha|\le 9}\|\Gamma^\alpha u'(s,\cd)\|_2 \le
C'A_0\varepsilon^2\langle s\rangle^{-1+A_2\varepsilon}.
\end{equation*}
Therefore
$$\int_0^t \sum_{|\alpha|\le 8}\|\Gamma^\alpha G(s,\cd)\|_2\, ds
\le C'A_0\varepsilon^2 \int_0^t\langle
s\rangle^{-1+A_2\varepsilon}\, ds \le
\frac{C'A_0}{A_2}\varepsilon(1+t)^{A_2\varepsilon}.$$ Since this
give the desired bounds for the remaining terms in \eqref{4.4} if
$A_2$ is large enough, the proof is complete. \qed

\end{document}